\documentclass[3p,times]{elsarticle}
\usepackage{amsmath}
\usepackage{enumitem}
\usepackage{pifont}
\usepackage[utf8]{inputenc}
\usepackage[colorlinks]{hyperref}

\makeatletter
\def\LaTeX{\leavevmode L\raise.42ex
    \hbox{\kern-.3em\size{\sf@size}{0pt}\selectfont A}\kern-.15em\TeX}
\makeatother

\newcommand{\BibTeX}{{\rm B\kern-.05em{\sc
          i\kern-.025emb}\kern-.08em\TeX}}
\makeatletter
\def\@currentlabel{2.1}\label{e:dispaa}
\def\@currentlabel{2.21}\label{e:dispau}
\def\@currentlabel{2.22}\label{e:dispav}
\def\@currentlabel{2.23}\label{e:dispaw}
\def\@currentlabel{2.24}\label{e:dispax}
\def\theequation{\thesection.\@arabic\c@equation}
\makeatother

\renewcommand{\theequation}{\arabic{section}.\arabic{equation}}

\newcommand{\R}{\mathbb R}

\newcommand{\N}{\mathbb N}

\def \D{\Delta}
\def \O{\Omega}
\newtheorem{theorem}{Theorem}

\newtheorem{theo}{Theorem} [section]
\newtheorem{lem}{Lemma} [section]

\newtheorem{definition}{Definition} [section]

\newtheorem{rem}{Remark}[section]

\renewcommand{\theequation}{\thesection.\arabic{equation}}
\renewcommand{\thesection}{\arabic{section}}
\renewcommand{\theequation}{\thesection.\arabic{equation}}
\let\ssection=\section\renewcommand{\section}{\setcounter{equation}{0}\ssection}
\begin{document}
\begin{frontmatter}
\title{Multiple solutions for Grushin operator without odd nonlinearity}
\author[mk,ah]{Mohamed Karim Hamdani \corref{cor1}}
\ead{hamdanikarim42@gmail.com}
\cortext[cor1]{Corresponding author:Tel:+21624124651.}
\address[mk]{Mathematics Department, University of Sfax, Faculty of Science of Sfax, Sfax, Tunisia.}
\address[ah]{Military School of Aeronautical Specialities, Sfax, Tunisia.}
\begin{abstract}
We deal with existence and multiplicity results for the following nonhomogeneous and homogeneous equations, respectively:
\begin{eqnarray*}
(P_g)\quad - \Delta_{\lambda} u + V(x) u = f(x,u)+g(x),\;\mbox{ in } \R^N,\;
\end{eqnarray*}
and
\begin{eqnarray*}
(P_0)\quad - \Delta_{\lambda} u + V(x) u = K(x)f(x,u),\;\mbox{ in }  \R^N,\;
\end{eqnarray*}
where  $\Delta_{\lambda}$  is the strongly degenerate
operator, $V(x)$ is allowed to be sign-changing, $K\in C(\R^N,\R)$, $g:\R^N\to\R$ is a perturbation and the nonlinearity $f(x,u)$  is
a continuous function does not satisfy the Ambrosetti-Rabinowitz superquadratic
condition ($(AR)$ for short). First, via the mountain pass
theorem and the
Ekeland’s variational principle, existence of two different solutions for $(P_g)$ are obtained when $f$ satisfies superlinear growth condition.  Moreover, we prove the existence of infinitely many solutions for $(P_0)$ if $f$ is odd in $u$ thanks an extension of Clark's theorem near the origin. So, our main results considerably improve results appearing in
the literature.
\end{abstract}
\end{frontmatter}

\section{Introduction}
In this paper we consider the following equations  in both nonhomogeneous case $g(x)\neq 0$, namely
\begin{eqnarray*}
(P_g)\quad - \Delta_{\lambda} u + V(x) u = f(x,u)+g(x),\;\mbox{ in }  \R^{N_1}\times \R^{N_2} = \R^N,\;
\end{eqnarray*}
and in the homogeneous case $g(x)\equiv 0$, that is
\begin{eqnarray*}
(P_0)\quad - \Delta_{\lambda} u + V(x) u = K(x)f(x,u),\;\mbox{ in }  \R^{N_1}\times \R^{N_2} = \R^N,\;
\end{eqnarray*}
where $N\geq 2$, $V:\R^N\to \R$ is a potential function which is allowed to be sign-changing, $K\in C(\R^N,\R)$, $g:\R^N\to\R$ is a perturbation and $f:\R^N\times \R \to\R$ is
a continuous function and $ \Delta_{\lambda}$ is the Grushin operator defined by $\Delta_{\lambda} u= \Delta_x u+|x|^{2\lambda}\Delta_y u$ with $\lambda> 0$ $(x\in \R^{N_1}, y\in \R^{N_2})$. Let $\nabla_\lambda$ be the gradient operator defined by $$\nabla_\lambda u=(\nabla_xu,|x|^\lambda\nabla_yu) ~\mbox{ and }~ |\nabla_\lambda u|^2=|\nabla_xu|^2+|x|^{2\lambda}|\nabla_yu|^2.$$
We have $$\Delta_{\lambda}=\nabla_\lambda\cdot\nabla_\lambda.$$

Recently, a lot of attention has been paid to the study of the existence and multiplicity results for problems $(P_g)$ (resp. $P_0$). For this reason, many authors have devoted their attention to problems of this types and they have widely studied them by using variational methods under various conditions on the potentials $V(x)$, $g(x)$ and $K(x)$ and the nonlinearity $f(x,u)$ especially when it is
 superlinear or sublinear at infinity. The paper was motivated by some works appeared in recent years. Let us start with
the nonhomogeneous case $g(x)\neq 0$. In \cite{CL}, Chen-Li studied the following problem in $\R^N$ \begin{eqnarray}\label{CL}
-\left(a+b\int_{\R^N}|\nabla u|^2dx\right)\Delta u +V(x) u=f(x,u)+g(x),\;\mbox{ in } \R^N,
 \end{eqnarray}
 where $a>0$, $b\geq 0$, $V(x)$ is allowed to be sign-changing, $g\in C^1(\R^N)\cap L^2(\R^N)$ is a radial function with $g\neq 0$ and $f$ satisfying the following variant Ambrosetti-Rabinowitz type condition ($(AR)$ for short):
 \\\\
{\bf(AR)}\;\;there exist $\mu>4$ such that
$$\mu F(x,s):=\mu\int_{0}^{s}f(x,t)dt \leq sf(x,s), \;\mbox{ for any } s\in\R \mbox{ and } x\in \R^N.$$
They showed that there are two constants $m_0>0$ (respectively $m_1 > 0$) such that problem \eqref{CL} admits at least two different solutions  when $\|g\|_{L^2(\R^N)} < m_0$ (respectively at least two different radial solutions in $H^1(\R^N)$ when $\|g\|_{L^2(\R^N)} < m_1$ if $V(x)\equiv 1$ and $f(x,u)=|u|^{p-2}u,$ with $4<p<2^*$).\footnote{We refer the readers to \cite{Ham} for a new kirchhoff problem with $(AR)$ condition.}

It is well known that under $(AR)$, any Palais Smale
sequence of the corresponding energy functional is bounded.  Although $(AR)$ condition is very crucial, there are many
super-linear functions which do not satisfy the $(AR)$ condition. For instance the
function
$$f(x,s)=s^3\log(1 + |s|)$$
does not satisfy the $(AR)$ condition. Because of this reason some attempts were made to replace condition $(AR)$. For example, in \cite{ZZ}, Zhang-Xu also studied  this
problem without the $(AR)$ condition, $a,b$ are positive constants, $V(x)\equiv 1$ and $f(x,u)=|u|^{p-2}u$ with $p\in(1,5)$, \footnote{note that when $p\in (1,3]$,  $|u|^{p-1}u$ do not satisfies $(AR)$ condition,} and $g\in C^1(\R^3)\cap L^2(\R^3)$ such that
$0\leq g(x)=g(|x|)\in L^2(\R^3) \mbox{ and } \|g\|_{L^2(\R^3)}<m_3,$ where $m_3$ is a positive constant. The authors investigated the existence of at least two nontrivial radial solutions which the first solution with
negative energy obtained by using Ekeland's variational principle and the
second solution has positive energy by using the mountain pass geometry structure.

Regarding the homogeneous case $g(x)\equiv 0$. In $2016$, Li-Zhong \cite{LZ} investigated the existence of infinitely many solution for the following equation
\begin{eqnarray}\label{sd}
-\left(1+b\int_{\R^3}(\nabla u|^2+V(x)u^2)dx\right)[\Delta u+V(x)u]=K(x)f(x,u)\;\mbox{ in }\R^3,
 \end{eqnarray}
where $b>0$ is a constant, $K(x)\equiv 1$, $V(x)$ is a nonnegative potential function and the nonlinear term $f(x,u)$ is only locally defined for $|u|$ small and satisfies some mild conditions.\footnote{Without any growth conditions imposed on $f(x,u)$ at infinity with respect to $u$.}

Few years later, in $2015$, Feng-Feng \cite{FF} considered a class of Kirchhoff type problem like
equation \eqref{sd} with $V(x)\equiv 1$, $K\in L^\frac{2}{2-\gamma}(\R^3)\cap L^\infty(\R^3)$ is a positive continuous function and $f(x,u)$ satisfes sublinear condition in $u$ at infinity. By using the iterating method in \cite{LW}, they proved the following estimate which is valid to get the infinitely of many solutions via Clark's theorem such that:
\begin{eqnarray}\label{estimate}
\|u\|_{L^\infty(\R^3)}\leq C_1\|u\|^v_{L^6(\R^3)},
\end{eqnarray}
where $v$ is a number in $(0,1)$ and $C_1 > 0$ is independent of $u$.

Motivated by the previously mentioned works, elliptic problems involving the Grushin operator have attracted the attention of many authors, see e.g. \cite{R2,RH,T0,T1} and references therein. See also \cite{AM, CTG, KL, TT1, TT2} for results related to a more general class of degenerate operators, namely the $\Delta_\lambda$-Laplace operators. In \cite{KL}, Kogoj and Lanconelli investigated the $\Delta_\lambda$-Laplace operator under the
additional assumption that the operator is homogeneous of degree two with respect to a semigroup of
dilations in $\R^N$. In $2015$, Anh and My studied in \cite{AM} the following problem
\begin{eqnarray}\label{Pp}\begin{cases} - \Delta_{\lambda} u +V(x) u= f(x, u)&\text{ in }  \Omega\\
u=0 &\text{ in } \partial\Omega,
\end{cases}
\end{eqnarray}
where $\O$ is a bounded domain in $\R^N$ $(N \geq 2)$, $\Delta_\lambda$ is defined below, $V\equiv 1$ and $f(x, u)$ is a Carathéodory function which satisfies some subcritical growth and does not require the usual $(AR)$ condition. By the mountain pass theorem, the authors showed that the considered
problem admits at least one nontrivial weak solution and by the Fountain
theorem, infinitely many weak solutions. In $2017$, Chen-Tang-Gao \cite{CTG} studied \eqref{Pp} with $V(x)$ is allowing to be sign-changing such that
$$V\in C(\O,\R)\mbox{ and } \inf_{x\in \O}V(x)>-\infty\mbox{ for all } x\in \O,$$
 and $f$ is a function with a more general super-quadratic growth, which is weaker than the $(AR)$ condition.
By symmetric mountain pass theorem, the authors proved the existence of infinitely many solutions for problem \eqref{Pp}. Inspired by \cite{CTG}, Rahal-Hamdani \cite{RH} studied the following problems
\begin{eqnarray}\label{P'1}\begin{cases} - \Delta_{\lambda} u + {V}(x) u = {f}(x, u)+g(x)|u|^{q-2}u &\text{ in }  \Omega\\
u=0 &\text{ on } \partial\Omega,
\end{cases}
\end{eqnarray}
and
\begin{eqnarray}\label{P'2}\begin{cases} - \Delta_{\lambda} u + {V}(x) u = {f}(x, u)  +\lambda u&\text{ in }  \Omega\\
u=0 &\text{ on } \partial\Omega,
\end{cases}
\end{eqnarray}
where $\lambda\in\R$, the potential ${V}$ and the nonlinearities ${f}$ and $g$ satisfy the following conditions:
\begin{itemize}
  \item[$(V1)$] $\displaystyle{\inf_{\O}} {V}(x)\geq V_1>0,$ and meas$\{x\in\O:\;{V}(x)\leq M\}<+\infty$, $\forall\;M>0$.
  \item [$(H_1)$] ${f} \in C(\O\times \R,\R)$ and  there exist $d_1 > 0$ and $d_2 > 0$ such that
      $$|{f}(x,t)|\leq d_1|t|+d_2|t|^{s-1},\forall(x,t)\in \O\times\R,\text{ where }s\in(2,2^*_\lambda).$$
  \item [$(H_2)$] There exists $\mu > 2$ such that
  $$0<\mu {F}(x,t)\leq{f}(x,t)t,\;\forall\; |t|\geq r,x\in\O,$$
  where $ {F}(x,t)= \displaystyle \int_0^t  {f}(x,s)ds.$
  \item [$(H_3)$] $ {f}(x,t)=o(|t|)$ as $t\rightarrow 0$ uniformly in $x\in\O.$
  \item [$(H_4)$] $ {f}(x,-t)=- {f}(x,t)$ for all $(x,t)\in\O\times \R.$
\item [$(g)$]  $0\leq g \in L^{q'}(\O)\;\mbox{ with } q'=\frac{s}{s-q},\; \mbox{ where } \; q\in(1,2),s\in(2,2^*_\lambda) \mbox{  and } 2^*_\lambda=\frac{2Q}{Q-2}$.
\end{itemize}
By using the Fountain theorem and the mountain pass theorem, the authors in \cite{RH} achieved the following theorems.
\begin{theo}
   Suppose that $(V1)$, $(H_1)-(H_3)$ and $(g)$ hold. Then
there exists $\lambda_0>0$ such that if $||g||_{L^{q'}} <\lambda_0$, the problem \eqref{P'1} has at least one nontrivial solution.
\end{theo}

\begin{theo}
Suppose that $(V1)$, $(H_1)$, $(H_2)$, $(H_4)$ and $(\Lambda)$ hold.
Then the problem \eqref{P'2} has infinitely many nontrivial solutions $\{u_k\}^\infty_{k=1}$ with unbounded
energy.
\end{theo}
\begin{rem}
Rahal-Hamdani \cite{RH}  studied \eqref{P'1} and \eqref{P'2} with sign-changing potential ${V}$ and with the Ambrosetti-Rabinowitz condition. They got only for the superlinear case the existence of at least one nontrivial solution and multiple solutions with odd nonlinearity. Here, for the nonlinearity of $(P_g)$ {\bf we do not need to be odd and does not satisfy the Ambrosetti–Rabinowitz condition} and then we obtain for the superlinear case {\bf two solutions} for problem $(P_g)$ and for the sublinear situation case we prove  that $(P_0)$ has infinitely many solutions.
\end{rem}

Motivated by the papers mentioned above, in the first part of this paper, we need to make the following assumption on the potential $V(x)$:
\begin{itemize}
\item[$(\widetilde{V_1})$] $\tilde{V}\in C(\R^N,\R)$ and $\inf_{\R^N}\tilde{V}(x)\geq V_1>0,$ and meas$\{x\in\R^N:\;\tilde{V}(x)\leq M\}<+\infty$, $\forall\;M>0$.
\end{itemize}

Similarly to \cite{CTG}, we choose a constant $V_0 >
0$ such that $\tilde{V
}(x) := V(x) + V_0 \geq V_1> 0$ for all $x \in \R^N$ and let $\tilde{f}(x, u) = f (x, u) + V_0u, \forall (x, u) \in \R^N \times \R$. Then we obtain the following problem which is equivalent to $(P_g)$
\begin{eqnarray}\label{P2}- \Delta_{\lambda} u + \tilde{V}(x) u= \tilde{f}(x, u)+g(x) &\text{in }  \R^N.
\end{eqnarray}
In what follows, we turn our aim to study the problem \eqref{P2}.

The first aim of this paper is to study $(P_g)$ under more generic superlinear growth conditions in order to generalize or to give complementary results to the ones listed above. More precisely, our methods are different
than the method used in \cite{CTG,RH}, and we investigate existence of two different solutions for $(P_g)$ which one is negative energy solution and the other is positive
energy solution under the following conditions on $\tilde{f}$ and $g$:
\begin{itemize}
  \item [$(H_1)$] $\tilde{f} \in C(\R^N\times \R,\R)$ and  there exist $c_1 > 0$ and $q\in (4,2^*_\lambda)$ such that
      $$|\tilde{f}(x,t)|\leq c_1(1+|t|^{q-1}),\forall(x,t)\in \R^N\times\R,$$
where $2^*_\lambda=\frac{2Q}{Q-2}$ and $Q$ denotes the homogeneous dimension of $\R^N$ with
respect to a group of dilations $\{\delta_t\}_{t>0},$ i.e.
$$Q:=\epsilon_1+...+\epsilon_N.$$
  \item [$(H_2)$] $\lim_{|t|\to +\infty}  \frac{\tilde{F}(x,t)}{|t|^4}=+\infty$,  uniformly for $x\in \R^N$,\\ where $\tilde{F}(x,t)=\int_{0}^{t}\tilde{f}(x,s)ds$.
    \item [$(H_3)$] $\tilde{f}(x,t)=o(|t|)$ as $|t|\rightarrow 0$ uniformly in $x\in\R^N;$
    \item[$(H_4)$]  there exists  $C>0$ such that $C|\tilde{f}(x,s)|^{(2_\lambda^*)'}\leq s \tilde{f}(x,s)- 2\tilde{F}(x,s)$,
 $\forall s\in \R$ and $x\in \R^N,$
 \\where $(2_\lambda^*)'=\frac{2_\lambda^*}{2_\lambda^*-1}$ the conjugate exponent of $2^*_\lambda$.
 \item [$(G)$]  $g \in L^2 (\R^N)$ and $g(x)\geq 0 (\not\equiv 0)$ for almost every $x\in \R^N.$
\end{itemize}
Now, we are ready to state our first main result
\begin{theo}\label{theorem 1.2}
  Suppose that $(\widetilde{V_1})$  and $(H_1)-(H_4)$ are satisfied and suppose that
$g\in L^{2}(\R^N)$. Then there exists a constant $\delta_0 > 0$ such that problem $(P_g)$ admits
at least two different nontrivial solutions in $E_g$, provided that $\|g\|_{L^{2}(\R^N)}\leq \delta_0.$
\end{theo}
\begin{rem}
We note that condition like $(H_4)$ was first
introduced by Harrabi in \cite{H} for the polyharmonic problems and by Hamdani-Harrabi \cite{HH}  for the $m$-polyharmonic Kirchhoff problems (see also Hamdani \cite{H2} for the asymmetric $m-$laplacian Kirchhoff problems, Harrabi \cite{H2} for the fourth-order elliptic Equations and Hajlaoui-Harrabi \cite{H3} for a priori estimates and existence of positive solutions for higher-order elliptic equations) which is weaker than $(H_2)$ of \cite{RH}.
\end{rem}

In the second part of the paper we study $(P_0)$ under a sublinear situation. \footnote{We strongly need
that perturbation  $g\equiv 0$ in the equation $(P_0)$ because this requirement is strictly necessary to solve the delicate lack of compactness.} Furthermore, we assume that the potential $V$ and the perturbation function $K$ satisfy the following:
\begin{enumerate}[label=($V\sb{\arabic*}$):,
ref=$V\sb{\arabic*}$]
  \item $V,K\in C(\R^N,\R), V(x)\geq \alpha_0$ and $0<K(x)\leq \beta_0$ for some $ \alpha_0>0,\beta_0>0$, and
      $M:=K^{\frac{2}{2-\gamma}}V^{\frac{-\gamma}{2-\gamma}}\in L^1 (\R^N)$\label{V1}
      for $\gamma\in [1,2)$.
\end{enumerate}
Next, we assume for $\delta>0$ that the nonlinearity $f:\R^N\times [-\delta,\delta]\to\R$ is a continuous function satisfying:
\begin{enumerate}[label=($f\sb{\arabic*}$):,
ref=$f\sb{\arabic*}$]
  \item $f$ is odd in $t$, i.e. $f(x,-t)=-f(x,t)$ for all $x\in \R^N$ and $t\in \R$;\label{f1}
  \item there exist $\gamma\in [1,2)$ and $C>0$ such that $|f(x,t)|\leq C|t|^{\gamma-1}$;\label{f2}
  \item $\lim_{t\to 0}\frac{F(x,t)}{|t|^2}=+\infty$ uniformly in some ball $B_r(x_0)\subset\R^N$,\\
       where $F(x,t)=\int_{0}^{t}f(x,s)ds$.\label{f3}
\end{enumerate}
Now, we state our second main results which generalizes the main results of \cite{FF}:
\begin{theo}\label{thm1}
Suppose that \eqref{f1}-\eqref{f3} and \eqref{V1} are satisfied. Then problem $(P_0)$ has infinitely
many solutions $\{u_k\}$ such that $\|u_k\|_{L^\infty}\to 0$ as $k\to \infty$.
\end{theo}

Finally, let us simply describe the main approaches to obtain Theorems \ref{theorem 1.2} and \ref{thm1}. To show the
existence of at least two different energy solutions in the superlinear case which one is
negative energy solution and the other is positive energy solution, we shall use the mountain pass theorem of Rabinowitz \cite{R}  and  Ekeland’s variational principle of  Mawhin-Willem \cite{MW} stated in Section \ref{Plusieurs}. Differently to  \cite{FF} in the sublinear case (see estimate \eqref{estimate}), here we will prove the following estimate
 \begin{eqnarray*}
\|u\|_{L^\infty(\R^N)}\leq C_1\|u\|_{L^{2_\lambda^*}(\R^N)}^{v},
\end{eqnarray*}
where $v=\prod_{i=0}^{\infty} \frac{\alpha_i+\gamma}{\alpha_i+2}$ is a number in $(0,1)$, $\alpha_i$ is a positive number and
$C_1=\exp\left( \sum_{i=0}^{k}\frac{2\log(\sigma(\alpha_i+2))}{\alpha_i+2} \right)$,
for some $\sigma\geq 1$  which derive the existence of infinitely many solutions as well as the structure of the set of critical points near the origin.

To our best knowledge, Theorems \ref{theorem 1.2} and \ref{thm1} are new even in the study for the Grushin problem with sign-changing potential in $\R^N$.

This paper is organized as follows. In Section \ref{section2}, we give some preliminaries and notation. Section \ref{section3} is
devoted to the proof of Theorem \ref{theorem 1.2}. In section \ref{section4}, we prove Theorem \ref{thm1}.
\section{Preliminaries and notation}\label{section2}
\subsection{\bf Function spaces and embedding theorem}

We recall the functional setting in \cite{CTG, KL, KS}. We consider the operator of the form
$$\Delta_\lambda:=\sum_{i=1}^{N}\partial_{x_i}(\lambda_i^2\partial_{x_i}),$$
where $\partial_{x_i}=\frac{\partial}{\partial_{x_i}}$, $i=1,...,N$. Here the functions $\lambda_i : \R^N \rightarrow \R$ are continuous, strictly
positive and of class $C^1$ outside the coordinate hyperplanes, i.e. $\lambda_i>0,i=1,...,N$ in $\R^N\setminus\prod$, where $\prod=\{(x_1,...,x_N)\in\R^N:\prod_{i=1}^{N}x_i=0\}.$ As in \cite{KL} we assume that $\lambda_i$
satisfy the following properties:
\begin{enumerate}
  \item $\lambda_1(x)\equiv 1, \lambda_i(x)=\lambda_i(x_1,...,x_{i-1}),i=1,...,N;$
  \item  for every $x\in\R^N$, $\lambda_i(x)=\lambda_i(x^*), i=1,...,N;$
  where $x^*=(|x_1|,...,|x_N|) \mbox{ if } x=(x_1,...,x_N);$
  \item there exists a constant $\rho\geq 0$ such that
  $$0\leq x_k\partial_{x_k}\lambda_i(x)\leq \rho \lambda_i(x)\;\;\forall k\in\{1,...,i-1\},i=2,...,N,$$
and for  every $x\in\R_+^N:=\{(x_1,...,x_N)\in \R^N:x_i\geq0\;\forall i=1,...,N\};$
  \item there exists a group of dilations $\{\delta_t\}_{t>0}$
  $$\delta_t:\R^N\rightarrow\R^N,\delta_t(x)=\delta_t(x_1,...,x_N)=(t^{\epsilon_1}x_1,...,t^{\epsilon_N}x_N),$$
  where $1\leq \epsilon_1\epsilon_2\leq...\leq\epsilon_N$, such that $\lambda_i$ is $\delta_t-$homogeneous of degree $\epsilon_i-1$,
  i.e.
$$
\lambda_i(\delta_t(x))=t^{\epsilon_i-1}\lambda(x), \forall x\in \R^N, t>0,i=1,...,N.$$
  This implies that the operator $\Delta_\lambda$ is $\delta_t-$homogeneous of degree two, i.e.
$$\Delta_\lambda(u(\delta_t(x)))=t^2(\Delta_\lambda u)(\delta_t(x)),\;\;\forall u\in C^\infty(\R^N).$$
\end{enumerate}
Now, we denote by $W_\lambda^{1,2}(\R^N)$ the closure of $C_0^1(\R^N)$ in the norm
$$||u||_{1,2}:=\Big(\int_{\R^N}|\nabla_\lambda u|^2 dx\Big)^{\frac{1}{2}}.$$
In view of the presence of the potential $\tilde{V}(x)$  (resp. $V(x)$ in \eqref{V1}), we consider the space
$$E_g:=\Big\{u\in W_\lambda^{1,2}(\R^N):\;\;\int_{\R^N}\tilde{V}(x)|u(x)|^2dx<\infty\Big\},$$
$$(\mbox{resp. }E_0:=\Big\{u\in W_\lambda^{1,2}(\R^N):\;\;\int_{\R^N}V(x)|u(x)|^2dx<\infty\Big\}),$$
equipped with the following inner
product $$(u,v)=\int_{\R^N}(\nabla_\lambda u\cdot\nabla_\lambda v + \tilde{V}(x) uv)dx,$$
$$(\mbox{resp. } (u,v)=\int_{\R^N}(\nabla_\lambda u\cdot\nabla_\lambda v + V(x) uv)dx).$$
Then $E_g$ ({resp.} $E_0$) is a Hilbert space endowed with the norm
\begin{eqnarray*}
\|u\|:=\Big(\int_{\R^N}|\nabla_\lambda u|^2+\int_{\R^N}\tilde{V}(x)|u(x)|^2dx\Big)^{\frac{1}{2}}
\simeq \Big(||u||^2_{{1,2}}+||\sqrt{\tilde{V}}\; u||^2_{L^2(\R^N)}\Big)^{\frac12},
\end{eqnarray*}
$$(\mbox{resp. } \|u\|:=\Big(\int_{\R^N}|\nabla_\lambda u|^2+\int_{\R^N}V(x)|u(x)|^2dx\Big)^{\frac{1}{2}}
\simeq \Big(||u||^2_{{1,2}}+||\sqrt{V}\; u||^2_{L^2(\R^N)}\Big)^{\frac12}).$$

\subsection{\bf Examples of $\D_{\lambda}$ operators (see \cite{AE,RH} for other examples)}
$\\$
\textbf{Example 1.} \rm
Let $\alpha, \beta$ and $\gamma$ be nonnegative real constants.
 We consider the operator
\[
\Delta_\lambda =\Delta_{x^{(1)}} + |x^{(1)}|^{2\alpha} \Delta_{x^{(2)}}
+ |x^{(1)}|^{2\beta} |x^{(2)}|^{2\gamma} \Delta_{x^{(3)}},
\]
where $\lambda= (\lambda^{(1)},\lambda^{(2)},\lambda^{(3)})$ with
\begin{gather*}
\lambda_j^{(1)} (x) \equiv 1,\quad j=1,\dots,N_1 \\
\lambda_j^{(2)} (x) = |x^{(1)}|^{\alpha}, \quad j=1,\dots,N_2, \\
\lambda_j^{(3)}(x) = |x^{(1)}|^{\beta}|x^{(2)}|^{\gamma}, \quad j=1,\dots,N_3.
\end{gather*}
The dilations become
\[
\delta_r\Big(x^{(1)},x^{(2)},x^{(3)}\Big)
=\Big( r x^{(1)}, r^{\alpha+1} x^{(2)}, r^{\beta + (\alpha +1)\gamma +1} x^{(3)}\Big).
\]
Similarly, for operators of the form
\begin{align*}
\Delta_\lambda
&= \Delta_{x^{(1)}} + |x^{(1)}|^{2\alpha_{1,1}} \Delta_{x^{(2)}}+
 |x^{(1)}|^{2\alpha_{2,1}}|x^{(2)}|^{2\alpha_{2,2}} \Delta_{x^{(3)}} +\dots \\
&\quad +
 \Big( \prod_{i=1}^{k-1} |x^{(i)}|^{2\alpha_{k-1,i}}\Big)\Delta_{x^{(k)}},
\end{align*}
where $\alpha_{i,j}\geq 0$, $i=1,\dots,k-1, j=1,\dots, i$, are real constants,
the group of dilations is given by
\[
\delta_r \Big(x^{(1)}, \dots,x^{(k)}\Big)
=\Big(r^{\varepsilon_1} x^{(1)},\dots, r^{\varepsilon_k} x^{(k)}\Big)
\]
with $\varepsilon_1 =1$ and $\varepsilon_j =1+\sum_{i=1}^{j-1}\alpha_{j-1,i} \varepsilon_{i}$,
for $i=2,\dots,k$.
In particular, if
$\alpha_{1,1}=\dots=\alpha_{k-1,k-1}=\alpha$,
\[
\delta_r \Big(x^{(1)}, \dots, x^{(k)}\Big)
= \Big( r x^{(1)}, r^{\alpha+1} x^{(2)},\dots, r^{(\alpha+1)^{k-1}} x^{(k)}\Big).
\]
\\
\textbf{Example 2.} Let $\lambda$ be a real positive constant and $k = 2$. We consider the Grushin-type
operator $$ \D_{\lambda}= \D_x + |x|^{2 \lambda } \D_y,$$
where $\lambda = (\lambda_1,\lambda_2)$ with $$\lambda_1(x)=1, \quad \lambda_2(x)= |x^{(1)}|^{\lambda },\quad x \in \R^{n_1} \times \R^{n_2}.$$
Our group of dilations is
$$\delta_t(x)= \delta_t(x^{(1)},x^{(2)})= (tx^{(1)}, t^{\lambda+1}x^{(2)}),$$
and the homogenous dimension with respect to $(\delta_t)_{t>0}$ is $Q = n_1 + (\lambda + 1)n_2$.

\begin{lem}\label{Lem2.1}(See \cite{AM, CTG, KL}.) Evidently, $E_g$ (resp. $E_0$) is continuously embedded into $W_\lambda^{1,2}(\R^N)$ and hence continuously
embedded into $L^p(\R^N)$ for $1\leq p \leq 2_\lambda^*$ and the embedding from $E_g$ (resp. $E_0$) into $L^p(\R^N)$ is compact if $p\in [1, 2_\lambda^*)$, where $2^*_\lambda=\frac{2Q}{Q-2}$.
Consequently, there exists $\tau_p > 0$ such that
\begin{eqnarray}\label{injection}
  ||u||_p\leq \tau_p||u||,\;\;\forall u\in E_g\; (\mbox{resp.} E_0),
\end{eqnarray}
where $||u||_p$ denotes the usual norm in $L^p(\R^N)$.
\end{lem}
Now, we define a functional $I_g$ of $(P_g)$ (resp. $I_0$ of $(P_0)$) by
\begin{eqnarray}\label{f31}
I_g(u)=\frac12||u||^2-\int_{\R^N}\tilde{F}(x,u)dx-\int_{\R^N} g(x)udx,\;\forall u\in E_g,\\
(\mbox{resp. } I_0(u)=\frac12\|u\|^2-\int_{\R^N}K(x)F(x,u)dx,\;\forall u\in E_0)\nonumber.
\end{eqnarray}
The energy functional $I_g : E_g \rightarrow \R$ (resp. $I_0 : E_0 \rightarrow \R$) is well defined
and of class $C^1$. Moreover, the derivative of $I_g$ (resp. $I_0$) is
\begin{eqnarray}\label{f32}
\langle I_g'(u),v\rangle=\int_{\R^N}\nabla_\lambda u\nabla_\lambda vdx +\int_{\R^N}\tilde{V}(x)uvdx-\int_{\R^N}\tilde{f}(x,u)vdx
-\int_{\R^N} g(x)vdx,\\
(\mbox{resp. } \langle I'_0 (u),v\rangle=\int_{\R^N}\nabla_\lambda u\nabla_\lambda vdx +\int_{\R^N}V(x)uvdx-\int_{\R^N}k(x)F(x,u)dx)\nonumber,
\end{eqnarray}

for all $u,v\in E_g$ (resp. in $E_0$). Therefore, the critical points of $I_g$ (resp. $I_0$) are weak solutions for $(P_g)$ (resp. $(P_0)$).
\subsection{\bf Ekeland's variational principle, Mountain Pass theorem and Clark's theorem}\label{Plusieurs}
\begin{definition}
A sequence $\{u_n\} \subset E$ is said to be a $(P S)$ sequence if
\begin{eqnarray}\label{eqt3.10}
I_g (u_n) \to c \mbox{ and } I_g'(u_n)\to0  \mbox{ as }
n\to \infty \;(\mbox{resp. } I_0 (u_n) \to c \mbox{ and } I_0'(u_n)\to0  \mbox{ as }
n\to \infty),
\end{eqnarray}
where $c \in\R$. $I_g$ (\mbox{resp. } $I_0$) is said to satisfy the $(P S)$ condition if any $(P S)$ sequence has a
convergent subsequence.
\end{definition}

In the superlinear case, we will use the following version of Ekeland's variational principle and the mountain pass theorem to prove the existence of two different solutions. One is
negative energy solution and the other is positive energy solution.
\begin{theorem}(\cite{MW}, Ekeland's variational principle) Let $X$ be a complete metric space with
metric $d$ and let $I: X\to (-\infty, +\infty]$ be a lower semicontinuous function, bounded from
below and not identical to $+\infty$. Let $\epsilon> 0$ be given and $u \in X$ be such that
$$I(u) \leq \inf_X I +\epsilon.$$
Then there exists $v\in X$ such that
$$I (v) \leq I (u),\;\; d(u, v)\leq 1,$$
and for each $w\neq v$ in $X$, one has
$$I (v)-\epsilon d(v,w)< I (w).$$
\end{theorem}
\begin{theorem}(\cite{R}, mountain pass theorem)\label{mountain} Let $X$ be a real Banach space and $I\in  C^1(X,\R)$ satisfying $(PS)$ condition. Suppose $I (0)= 0$ and
\begin{enumerate}
  \item there exist two constants $\beta_0,\alpha_0>0$ such that $I \mid\partial Q_{\beta_0}\geq\alpha_0$;
  \item there is $u_1 \in X\mid \tilde{Q}_{\beta_0}$ such that $I(u_1) \leq0.$
\end{enumerate}
Then, $I$ possesses a critical value $c\geq \alpha_0$. Moreover, $c$ can be characterized as
$$c=\inf_{\chi \in \Gamma} \max_{u\in\chi([0,1])}I(u),$$
where
$\Gamma=\{\chi\in C([0,1],X):\chi(0)=0, \;\chi(1)=u_1\}.$
\end{theorem}
In the sublinear case, we give the improved Clark' theorem in \cite{LW} to prove Theorem \ref{thm1}.
\begin{theorem}(\cite{LW}  Clark theorem) \label{clark} Let $X$ be a real Banach space; $I \in C^1(X, \R)$ satisfies the $(PS)_c$ condition is even and bounded from below, and $f(0)=0$. If for any $k\in \N$, there
exists a $k-$dimensional subspace $X^k$ of $X$ and $\rho_k > 0$ such that $\sup_{X^k\cap S_\rho}I < 0$, where
$S_\rho
= \{u \in X  \mbox{ such that } \|u\|=\rho\}$, then at least one of the following conclusions holds:
\begin{enumerate}
  \item there exists a sequence of critical points $u_k$ satisfying $I(u_k) < 0$ for all $k$ and $\|u_k\|\to 0$
as $k \to \infty$;
  \item there exists $r > 0$ such that for any $0 < a < r$ there exists a critical point u such that
$\|u\|=a$ and $I (u) = 0$.
\end{enumerate}
\end{theorem}
\section{The superlinear case}\label{section3}
\begin{lem}\label{lemma32}
Assume that $(\widetilde{V_1})$, $( H_1)$ and $(H_3)$ hold. Then there exist some constants
$\beta_0$, $\alpha_0$, $\delta_0 > 0$ such that $I_g (u)\geq \alpha_0$ whenever $\|u\|\geq \beta_0$ and all $g\in L^2(\R^N )$, with $||g||_{L^2(\R^N)}<\delta_0$.
\end{lem}
{\bf Proof:} By $(H_1)$ and $(H_3)$, for all $\epsilon > 0$, there exists $C_\epsilon > 0$ such that

\begin{eqnarray}\label{f311}
|\tilde{f}(x,u)|\leq \epsilon|u|+C(\epsilon)|u|^{q-1}, \;\;\;\forall(x,u)\in\R^N\times\R,
\end{eqnarray}
and thus
\begin{eqnarray}\label{e312}
|\tilde{F}(x,u)|\leq \frac \epsilon 2|u|^2+\frac{C(\epsilon)}{q}|u|^{q}, \;\;\;\forall(x,u)\in\R^N\times\R.
\end{eqnarray}
From the Hölder inequality, \eqref{e312} and \eqref{injection}, we have
\begin{eqnarray*}
  I_g(u) &=& \frac{1}{2} ||u||^2-\int_{\R^N}\tilde{F}(x,u)dx- \int_{\R^N} g(x)udx\\
  &\geq&\frac 1 2 \|u\|^2-\frac \epsilon 2\|u\|^2_{L^2(\R^N)}-\frac{C(\epsilon)}{q}\|u\|_{L^q(\R^N)}^{q}- \|g\|_{L^2(\R^N)}\|u\|_{L^2(\R^N)}\\
  &\geq&\frac 12 ||u||^2-\frac\epsilon 2 \tau_2^2\|u\|^2-\frac{C(\epsilon)}{q}\tau_q^q||u||^q- \tau_2\|g\|_{L^2(\R^N)}\|u\|\\
  &\geq&\|u\|\left[\left(\frac12-\frac\epsilon 2 \tau_2^2\right)\|u\|-\frac{C(\epsilon)}{q}\tau_q^q\|u\|^{q-1}- \tau_2\|g\|_{L^2(\R^N)}\right].
\end{eqnarray*}
Taking $\epsilon=\frac{1}{2\tau_2^2}$ and let
$$h(t)=\frac14 t-\frac{C(\epsilon)}{q}\tau_q^qt^{q-1}\;\mbox{ for all } t\geq 0.$$
Note that $4<q<2^*_\lambda$, we can conclude that there exists a constant $\beta_0 > 0$ such
$$h(\beta_0)=\max_{t\geq 0}h(t)>0.$$
Taking $\delta_0:=\frac{h(\beta_0)}{\tau_2}$, we can get
$$I_g(u)\geq \frac{\beta_0h(\beta_0)}{2}:=\alpha_0>0,$$
where $\|u\|=\beta_0$ and $\|g\|_{L^2}\leq \delta_0$. This completes the proof.\qed
\begin{lem}\label{Lemma33}
  Assume that $(\widetilde{V_1})$ and $(H_1)–(H_3)$ hold. Then there exists a function
$w \in E_g$ with $\|w\|> \beta_0$ such that $I_g(w) < 0$.\end{lem}
{\bf Proof.}
By $(H_1)–(H_3)$, for any $M > 0$, there exists $C(M) > 0$ such that
\begin{eqnarray}\label{e38}
\tilde{F}(x,t)\geq M|t|^4-C(M)|t|^2,\;\forall(x,t)\in \R^N\times \R.
\end{eqnarray}
As $g \in L^2 (\R^N)$ and $g\geq 0 (\not\equiv 0)$, we can choose a function $\phi\in E$ such that
\begin{eqnarray}\label{ghgh}
\int_{\R^N} g(x)\phi dx>0.
\end{eqnarray}
Hence, from \eqref{e38}, \eqref{ghgh} and H\"older inequality, we can
get as $t\to \infty$
\begin{eqnarray*}
  I_g(t\phi) &=& \frac{t^2}{2} \|\phi\|^2-\int_{\R^N}\tilde{F}(x,t\phi)dx-t \int_{\R^N} g(x)\phi dx\\
  &\leq& \frac{t^2}{2} \|\phi\|^2-M t^4 \|\phi\|^4_{L^4(\R^N)}+C(M)t^2\|\phi\|^2_{L^2(\R^N)}\\
  &\leq& \frac{t^2}{2}-M t^4 \|\phi\|^4_{L^4(\R^N)}+C(M)t^2\|\phi\|^2_{L^2(\R^N)}\\
  &\to& -\infty.
\end{eqnarray*}
 Hence there exists $w = t_0\phi$ with $t_0 > 0$ large enough
such that $\|w\|  > \beta_0$ and $I_g(w) < 0$. This completes the proof.\qed

\begin{lem}\label{Lemma34} Let $(H_1)–(H_3)$ hold. Then $I_g$ satisfies the $(PS)$ condition.
\end{lem}
{\bf Proof.} We proceed by steps.\\
\textbf{Step 1.} We shall show that $u_n$ is bounded in $E_g$. First, from \eqref{f31}, we obtain
\begin{eqnarray*}
\int_{\R^N} g(x) u_ndx+\|u_n\|^2 =\langle I_g'(u_n),u_n\rangle+\int_{\R^N} f(x,u_n)u_n dx.
\end{eqnarray*}
From \eqref{injection} and \eqref{ghgh} and applying H\"older's inequality to the second term in the right-hand side and using , we obtain
\begin{eqnarray}\label{kirch hk}
\|u_n\|^2 \leq\int_{\R^N} g(x) u_ndx +\|u_n\|^2\leq \langle I_g'(u_n),u_n\rangle+ C\tau_2\left(\int_{\R^N} |f(x,u_n)|^{(2_\lambda^*)'}dx\right)^\frac{1}{(2_\lambda^*)'}\|u_n\|.\end{eqnarray}
In view of \eqref{f31}, \eqref{f32}  and $(H_4)$, it follows
\begin{align}\label{ing2}
\int_{\R^N}g(x) u_ndx+2I_g(u_n)-\langle I_g'(u_n), {u_n}\rangle&= \int_{\R^N}\left[ f(x,u_n)u_n -2F(x,u_n)\right]dx\nonumber\\&\geq C\int_{\R^N} |f(x,u_n)|^{(2_\lambda^*)'}dx.
 \end{align}
As a consequence of \eqref{injection} and \eqref{eqt3.10}, we also have
\begin{eqnarray*}\label{3.16}
\int_{\R^N} g(x) u_ndx+2I_g(u_n)-\langle I_g'(u_n), {u_n}\rangle\leq k(1+\|u_n\|),
\end{eqnarray*}
so that, by \eqref{ing2}, we have
\begin{eqnarray}\label{eelem4}
k(1+\|u_n\|)\geq\int_{\R^N} |f(x,u_n)|^{(2_\lambda^*)'}dx.
\end{eqnarray}
From \eqref{kirch hk}-\eqref{eelem4}, we arrive at the conclusion
\begin{equation*}
\|u_{n}\|^{2} \leq C( 1 + || u_{n} ||^{{\frac{1}{(2_\lambda^*)'}+1}} ).
\end{equation*}
 As $2> \frac{1}{(2_\lambda^*)'}+1$, then the $(PS)$ sequence $u_n$ is bounded in $E_g$.\\\\
\textbf{Step 2.} Here, we will prove that $\{u_n\}$ has
a convergent subsequence in $E_g$. It follows from Lemma \ref{Lem2.1}
that the embedding $$E_g\hookrightarrow L^p(\R^N)$$ is compact, where $1\leq p<2^*_\lambda$. Going if necessary to a subsequence, there exists $u\in E_g$ such that
\begin{eqnarray}\label{cvg}
u_n\rightharpoonup u\mbox{ in } E_g,\; u_n \to u \mbox{ in } L^p(\R^N),\;\; u_n(x)\to u(x), \mbox{ a.e. in } \R^N.
\end{eqnarray}
From \eqref{f311} and Hölder's inequality, we have
\begin{eqnarray*}
\Big|&\displaystyle\int_{\R^N}&(\tilde{f}(x,u_n)- \tilde{f}(x,u))(u_n-u)dx\Big|\\ &\leq& \int_{\R^N}|(\tilde{f}(x,u_n)- \tilde{f}(x,u))(u_n-u)|dx\\ &\leq & \int_{\R^N}[\epsilon|u_n|+\epsilon|u|+C(\epsilon)|u_n|^{q-1}
+C(\epsilon)|u|^{q-1}]|u_n-u|dx\\
&\leq& \epsilon\Big(\|u_n\|_{L^2(\R^N)}+\|u\|_{L^2(\R^N)}\Big)\|u_n-u\|_{L^2(\R^N)}
+C(\epsilon)\Big(\|u_n\|^{q-1}_{L^q(\R^N)}+\|u\|^{q-1}_{L^q(\R^N)}\Big)\|u_n-u\|_{L^q(\R^N)},
\end{eqnarray*}
which shows that
\begin{eqnarray}\label{cvgf}
\lim_{n\to \infty}\int_{\R^N}(\tilde{f}(x,u_n)- \tilde{f}(x,u))(u_n-u)dx=0.\end{eqnarray}

Now, by \eqref{cvg} and $(\widetilde{V_1})$, we have

\begin{eqnarray}\label{39v}
|\int_{\R^N} \tilde{V}(x)(u_n-u)^2dx|\leq \int_{\R^N} \tilde{V}(x)|(u_n-u)|^2dx\leq V_1\int_{\R^N} |(u_n-u)|^2dx\rightarrow 0 \mbox{ as } n\rightarrow\infty.
\end{eqnarray}
On the other hand, from the H\"{o}lder inequality, we have
\begin{align*}
\int_{\R^N}\tilde{V}(x)(u_n-u)^2dx&=\int_{\R^N}\tilde{V}(x)|u_n|^2dx-\int_{\R^N}\tilde{V}(x)u_nudx-\int_{\R^N}\tilde{V}(x)uu_ndx+\int_{\R^N}\tilde{V}(x)|u|^2dx\\
&\geq||\tilde{V}(x)^\frac 1 2u_n||_2^2-||\tilde{V}(x)^\frac12u_n||_2||\tilde{V}(x)^{\frac12}u||_2-||\tilde{V}(x)^\frac12u||_2\|\tilde{V}(x)^\frac12u_n\|_2
+\|\tilde{V}(x)^\frac12 u\|^2_2\\
&=\left(\|\tilde{V}(x)^\frac 1 2u_n\|_2-\|\tilde{V}(x)^\frac 1 2u\|_2\right)^2\geq 0.
\end{align*}
Hence, \eqref{39v} implies that
\begin{eqnarray}\label{cvgVV}
\|\tilde{V}(x)^\frac 1 2u_n\|_2\rightarrow \|\tilde{V}(x)^\frac 1 2u\|_2, \mbox{ as } n\rightarrow \infty.
\end{eqnarray}
Obviously, $\langle I_g'(u_n)-I_g'(u), u_n-u\rangle\to 0$ as $n \to \infty$, since $u_n\rightharpoonup u$ in $E_g$ and $I_g'(u_n)\to 0$
in $E_g^*$. Hence, \eqref{cvg}, \eqref{cvgf} and \eqref{39v} give as $n\to \infty$
\begin{align}\label{normmm}
o(1)&=\langle I_g'(u_n)-I_g'(u), u_n-u\rangle\nonumber\\&=\int_{\R^N}|\nabla_\lambda u_n- \nabla_\lambda u|^2dx +\int_{\R^N}\tilde{V}(x)\big(u_n-u\big)^2dx-\int_{\R^N}\big(\tilde{f}(x,u_n)-\tilde{f}(x,u)\big)\big(u_n-u\big)dx
\nonumber\\&=\int_{\R^N}|\nabla_\lambda u_n- \nabla_\lambda u|^2dx+o(1).
\end{align}
Therefore, from \eqref{cvgVV} and  \eqref{normmm} we have $\|u_n-u\|_E\to 0$ as $n\to \infty$. Since $E_g$ is a reflexive Banach space,
weak convergence and norm convergence imply strong convergence. Therefore, $u_n \rightarrow u$
strongly in $E_g$. This completes the proof of step $2$.

In conclusion, $I_g$ satisfies the $(PS)$ condition, as stated.\qed

\subsection*{\bf Proof of Theorem \ref{theorem 1.2}.}
The proof of this theorem is divided into two steps.
\\
{\bf Step 1: Existence of negative energy solution.}

 We will prove that there exists a function $u_0\in E_g$ such that $I_g'(u_0)=0$ and $I_g (u_0) < 0$.
\\
 By the proof of Lemma \ref{Lemma33} and $(H_1)–(H_3)$, there exist two constants $C_1, C_2 > 0$ such
that
\begin{eqnarray}\label{nnb}
\tilde{F}(x,t)\geq C_1|t|^4-C_2|t|^2, \; \forall (x,t)\in \R^N\times \R.
\end{eqnarray}
Hence, we obtain from $(G)$ and \eqref{nnb} that
\begin{eqnarray*}
  I_g(tv) &=&\frac{t^2}{2} \|v\|^2-\int_{\R^N}\tilde{F}(x,tv)dx- \int_{\R^N} g(x)v(x)dx \\&\leq& \frac{t^2}{2} \|v\|^2-C_1t^4\|v\|_{L^4(\R^N)}^4-C_2t^2\|v\|_{L^2(\R^N)}^2-t \int_{\R^N} g(x)v(x)dx<0
\end{eqnarray*}
for $t \in (0, 1)$ small enough, where $\rho_0 > 0$ is given in Lemma \ref{lemma32}. Thus, we get
$$c_0 = \inf\{I_g(u) : u \in \overline{B}_{\rho_0}
\} < 0,$$
where $B_{\rho_0}= \{u \in E_g : \|u\| < \rho_0\}$. By the Ekeland variational principle (see \cite{MW}) and
Lemma \ref{lemma32}, there exists a sequence $\{u_n\}_n \subset B_{\rho_0}$ such that

$$c_0 \leq I_g(u_n)\leq c_0 + \frac 1n
\mbox{ and } I_g(v) \geq I_g (u_n) -\frac 1n
\|v - u_n\|$$
for all $v \in \overline{B}_{\rho_0}$. Then a standard procedure gives that $\{u_n\}_n$ is a bounded (PS) sequence of $I_g$ .
Therefore, Lemmas \ref{lemma32} and \ref{Lemma34} imply that there exists a function $u_0\in B_{\rho_0}$ such that
\begin{eqnarray}\label{mlk}
 I_g'(u_0)=0
\mbox{ and } I_g(u_0)= c_0 < 0.
\end{eqnarray}
{\bf Step 2: Existence of positive energy solution.}
\\\\
 Next, we want to apply Theorem \ref{mountain} to prove the existence of the other solution.
 \\
  By Lemmas
\eqref{lemma32}-\eqref{Lemma33}, we know that all conditions of  the mountain pass theorem (see \cite{R}) are satisfied. Thus applying Lemma \ref{Lemma34},
we can conclude that there exists a function $u_1 \in E_g$ such that
\begin{eqnarray}\label{e330}
I_g'(u_1)=0 \mbox{ and } I_g (u_1)\geq \alpha_0> 0,
\end{eqnarray}
that is, $u_0$ is a positive energy solution. Therefore, it follows from \eqref{mlk} and \eqref{e330} that
$u_0
\neq u_1$. This completes the proof.\qed
\section{The sublinear case}\label{section4}
 In this section, we are ready to prove the Theorem \ref{thm1}. In the sequel,
for the sake of clarity, we divide the proof into several steps.
\\\\
{\bf Proof of step 1.}  Let $f_1(x,t) \in C(\R^N \times \R,\R)$ odd in $t \in \R$ be a functional example, so that
\begin{eqnarray}\label{definitionf1}
f_1(x, t)=
\begin{cases}
f(x, t), & \mbox{if }  x\in \R^N \mbox{ and } |t|<\frac{\gamma}{2}, \\
\frac{f(x,\frac{\gamma}{2})}{-\frac{\gamma}{2}}(t-\gamma), &\mbox{if }  x\in \R^N \mbox{ and } \frac{\gamma}{2}<t<\gamma \\ \frac{f(x,\frac{\gamma}{2})}{-\frac{\gamma}{2}}(t+\gamma), &\mbox{if }x\in \R^N \mbox{ and } -\gamma<t<-\frac{\gamma}{2} \\
  0, & \mbox{if } x\in \R^N \mbox{ and } |t|>\gamma.
\end{cases}
\end{eqnarray}
In order to investigate the existence of infinitely
many solutions for $(P_0)$, we shall apply Theorem \ref{thm1} to the above function $f_1(x, t)$ and its associated functional
\begin{eqnarray}\label{fun.energ2}
I_{0,1}(u)=\frac{\|u\|^2}{2}-\int_{\R^N}K(x)F_1(x,u)dx,
\end{eqnarray}
where $F_1(x,u)=\int_{0}^{u}f_1(x,s)ds$.
\\
From \eqref{f1} and \eqref{definitionf1}, it is easy to show that $I_{0,1}(u) \in C^1(E_0, \R)$, $I_{0,1}(u)$ is even, and $I_{0,1}(0) = 0$.

For all $u \in E_0$, we have
\begin{eqnarray*}
  \int_{\R^N}K(x)|F_1(x,u)|dx &\leq & C_1\int_{\R^N}K(x)|u|^\gamma dx  \\
 &\leq& C_1\|M\|_{L^1(\R^N)}^\frac{2-\gamma}{2}\|V|u|^2\|_{L^1(\R^N)}^\frac{\gamma}{2} \\
 &\leq& C_2\|u\|^\gamma.
\end{eqnarray*}
Therefore,
\begin{eqnarray*}
I_{0,1}(u)\geq \frac12\|u\|^2-C_2\|u\|^\gamma,\;\forall u\in E_0.
\end{eqnarray*}
Consequently, $I_{0,1}(u)$ is coercive and bounded from below.
\\\\
{\bf Proof of step $2$.} We claim that the functional $I_{0,1}(u)$ satisfies $(PS)$ condition in $E_0$.

Let $\{u_n\}\subset E_0$ be a $(PS)$ sequence for $I_{0,1}(u)$, that is
$$I_{0,1}(u_n)\to c,\quad I'_{0,1}(u)\to 0 \in E_0^*.$$

Then $\{u_n\}$ is bounded. Assume without loss of generality that {un} converges to $u$ weakly in
$E_0$, and by Lemma \ref{Lem2.1}, we may assume that
\begin{eqnarray}\label{cvg sub}
\begin{cases}
  u_n(x)\to u(x), & \mbox{a.e. in } B_R(0). \\\\
  u_n\to u, & \mbox{in } L^p( B_R(0)).
\end{cases}
\end{eqnarray}
where $p\in [1,2_\lambda^*).$
\\
For any $R > 0$, we have
\begin{eqnarray*}
\int_{\R^N} K(x)|f_1(x,u_n)&-&f_1(x,u)\|u_n-u|dx\\&\leq& c \int_{\R^N\backslash B_R(0)}K(x)(|u_n|^\gamma+|u|^\gamma)dx+c \int_{B_R(0)}(|u_n|^{\gamma-1}+|u|^{\gamma-1})(u_n-u) dx\\
&\leq& c\left(\|V|u_n|^2\|^\frac\gamma 2_{L^1({\R^N\backslash B_R(0)}) }+ \|V|u|^2\|^\frac\gamma 2_{L^1({\R^N\backslash B_R(0)}) }   \right)\|M\|^\frac{2-\gamma}{2}_{L^1({\R^N\backslash B_R(0)})}\\&+&c\left(\|u_n\|^{\gamma-1}_{L^\gamma(B_R(0))}+\|u\|^{\gamma-1}_{L^\gamma(B_R(0))}\right)\|u_n-u\|_{L^\gamma(B_R(0))},
\end{eqnarray*}
which implies that
\begin{eqnarray}\label{ee8}
\lim_{n\to \infty} \int_{\R^N} K(x)|f_1(x,u_n)-f_1(x,u)\|u_n-u|dx=0.
\end{eqnarray}
Now, by \eqref{cvg sub} and \eqref{V1}, we have
\begin{eqnarray}\label{38sub}
|\int_{\R^N} V(x)(u_n-u)^2dx|\leq \int_{\R^N} V(x)|(u_n-u)|^2dx\leq V_1\int_{\R^N} |(u_n-u)|^2dx\rightarrow 0 \mbox{ as } n\rightarrow\infty.
\end{eqnarray}
\\\\
On the other hand, from the H\"{o}lder inequality, we have
\begin{eqnarray*}
\int_{\R^N}V(x)(u_n-u)^2dx&=&\int_{\R^N}V(x)|u_n|^2dx-\int_{\R^N}V(x)u_nudx-\int_{\R^N}V(x)uu_ndx+\int_{\R^N}V(x)|u|^2dx\\
&\geq&\|V(x)^\frac 1 2u_n\|_2^2-\|V(x)^\frac12u_n\|_2\|V(x)^{\frac12}u\|_2-\|V(x)^\frac12u\|_2V(x)^\frac12u_n\|_2
+\|V(x)^\frac12 u\|^2_2\\
&=&\left(\|V(x)^\frac 1 2u_n\|_2-\|V(x)^\frac 1 2u\|_2\right)^2\geq 0.
\end{eqnarray*}
Hence, \eqref{38sub} implies that
\begin{eqnarray}\label{cvgV}
\|V(x)^\frac 1 2u_n\|_2\rightarrow \|V(x)^\frac 1 2u\|_2, \mbox{ as } n\rightarrow \infty.
\end{eqnarray}
Obviously, $\langle I_{0,1}'(u_n)-I_{0,1}'(u), u_n-u\rangle\to 0$ as $n \to \infty$, since $u_n\rightharpoonup u$ in $E_0$ and $I_{0,1}'(u_n)\to 0$
in $E_0^*$. Hence, \eqref{cvg sub}, \eqref{ee8} and \eqref{cvgV} give as $n\to \infty$
\begin{align}\label{norm}
o(1)&=\langle I_{0,1}'(u_n)-I_{0,1}'(u), u_n-u\rangle\nonumber\\&=\int_\O|\nabla_\lambda u_n- \nabla_\lambda u|^2 dx+\int_\O V(x)\big(u_n-u\big)^2dx-\int_{\R^N} K(x)(f_1(x,u_n)-f_1(x,u))(u_n-u)dx
\nonumber\\&=\int_\O|\nabla_\lambda u_n- \nabla_\lambda u|^2dx+o(1).
\end{align}
Therefore, from \eqref{cvgV} and  \eqref{norm} we have $\|u_n-u\|\to 0$ as $n\to \infty$. Since $E_0$ is a reflexive Banach space,
weak convergence and norm convergence imply strong convergence. Therefore, $u_n \rightarrow u$
strongly in $E_0$. This completes the proof of step $2$.
\\
In conclusion, $I_{0,1}$ satisfies the $(PS)$ condition, as stated.\qed
\\\\
{\bf Proof of step 3}

We show that problem $(P_0)$ has infinitely many Clark type solutions.
\\
By \eqref{f3}, we have that for any $K > 0$, there exists $\delta=\delta(K)> 0$ such that if $u\in C_0^\infty(B_r(x_0))$ and $|u|_\infty < \delta$ then $F_1(x,u)\geq K|u(x)|^2$ , and thus
$$I_{0,1}(u)=\frac{\|u\|^2}{2}-\int_{\R^N}K(x)|F_1(x,u)|dx\leq \frac12 \|u\|^2-K\|u\|_{L^2(\R^N)}^2.$$
This implies, for any $k\in \N$, if $X^k$ is a $k-$dimensional subspace of $C_0^\infty(B_r(x_0))$ and
$\rho_k > 0$ is sufficiently small then for any $u\in X^k\cap S_{\rho_k}$, there is a constant $C_k > 0$ such that $C_k\|u\|= C_k\rho_k < |u|_\infty <\delta$, where $S_{\rho_k}=\{u \in  X \mbox{ such that } \|u\|={\rho_k}\}$. This implies, for any $u\in X^k\cap S_{\rho_k}$
$$I_{0,1}(u)\leq \frac12 \|u\|^2-K\|u\|_{L^2(\R^N)}^2\leq (\frac12-K)\rho_k^2<0.$$
Now we
apply Theorem \ref{clark} to obtain infinitely many solutions $\{u_k\}$ for $(P_0)$ such that
\begin{eqnarray}\label{eqt23}
\|u_k\|\to0,\;k\to \infty.
\end{eqnarray}
Finally we show that $\|u_k\|_\infty\to0$ as $k\to \infty$.
Let $u$
a solution of $(P_0)$ and $\alpha>0$.
Let $T>0$  and set
 $u^T(x)=\max\{-T, \min\{u(x), T\}\}$.
 Multiplying both sides of $(P_0)$ with $|u^T|^\alpha u^T(x)$ implies
\begin{eqnarray*}
\int_{\R^N}-\Delta_\lambda u |u^T|^\alpha u^T(x) dx+\int_{\R^N}V u|u^T|^\alpha u^T(x)dx = \int_{\R^N}K(x) f_1(x,u)|u^T|^\alpha u^T(x)dx
\end{eqnarray*}
By the definition of $u^T$, we have that
\begin{eqnarray*}
\int_{\R^N}-\Delta_\lambda u |u^T|^\alpha u^T(x)dx&\geq& \int_{\R^N}-\Delta_\lambda u^T |u^T|^\alpha u^T(x)dx\\&=&
(\alpha+1)\int_{\R^N}\nabla_\lambda u^T(|u^T|^\alpha\nabla_\lambda u^T(x))dx\\&=&
(\alpha+1)\int_{\R^N}\left( \nabla_\lambda u^T(x)|u^T|^\frac\alpha 2\right)^2dx\\&=&
\frac{4(\alpha+1)}{(\alpha+2)^2}\int_{\R^N}\left| \nabla_\lambda \left| u^T\right|^{\frac{\alpha}{ 2}+1}\right|^2dx.
\end{eqnarray*}
From \eqref{f2} and \eqref{definitionf1}, we have
\begin{eqnarray*}
\frac{4(\alpha+1)}{(\alpha+2)^2}\int_{\R^N}\left| \nabla_\lambda \left| u^T\right|^{\frac{\alpha}{ 2}+1}\right|^2dx\leq \int_{\R^N}|u^T(x)|^{\alpha+\gamma}dx.
\end{eqnarray*}
By Lemma \ref{injection}, we get
\begin{eqnarray*}
  \|u^T\|_{L^\frac{(\alpha+2)N}{N-2}(\R^N)}\leq (\sigma(\alpha+2))^\frac{2}{\alpha+2} \|u^T\|_{L^{\alpha+\gamma}(\R^N)}^\frac{(\alpha+\gamma)}{(\alpha+2)},
\end{eqnarray*}
for some $\sigma\geq 1$ independent of $u$ and $\alpha$. Taking $\alpha_0=2_\lambda^*-1=\frac{Q+2}{Q-2}$ and $\alpha_k=\frac{(\alpha_{k-1}+2)Q}{Q-2}-1$, then $\alpha_k=\frac{(\frac{2_\lambda^*}{2})^{k+1}-1}{\frac{2_\lambda^*}{2}-1}\alpha_0$, for $k=1,2,....$
and in view
of the last inequality, an iterating process as in \cite{LW}
leads to
\begin{eqnarray*}
\|u^T\|_{L^{\alpha_{k+1}+1}(\R^N)}\leq \exp\left( \sum_{i=0}^{k}\frac{2\log(\sigma(\alpha_i+2))}{\alpha_i+2} \right)\|u^T\|_{L^{2_\lambda^*}(\R^N)}^{v_k},
\end{eqnarray*}
where $v_k=\prod_{i=0}^{k} \frac{\alpha_i+\gamma}{\alpha_i+2}.$ Sending $T$ to infinity then $k$ to infinity, consequently,
we derive
\begin{eqnarray*}
\|u\|_{L^\infty(\R^N)}\leq \exp\left( \sum_{i=0}^{k}\frac{2\log(\sigma(\alpha_i+2))}{\alpha_i+2} \right)\|u\|_{L^{2_\lambda^*}(\R^N)}^{v},
\end{eqnarray*}
where $v=\prod_{i=0}^{\infty} \frac{\alpha_i+\gamma}{\alpha_i+2}$ is a number in $(0,1)$ and
$\exp\left( \sum_{i=0}^{k}\frac{2\log(\sigma(\alpha_i+2))}{\alpha_i+2} \right)$ is a positive number. Therefore, $\|u_k\|_{L^\infty(\R^N)}\to 0$ as $k\to \infty$, and $u_k$ with $k$ sufficiently large
are solutions of $(P_0)$.\qed

\subsection*{\bf Acknowledgments}
The author would like to express his deepest gratitude to the Military School of Aeronautical Specialities, Sfax (ESA)
 for providing us with an excellent atmosphere for doing this work. Hamdani also wish to express his sincere gratitude to the anonymous referees for their knowledgeable report which provided insights that helped to improve the paper.
\center{
\end{document}